\documentclass{article}%
\usepackage{amssymb}
\usepackage{amsfonts}
\usepackage{graphicx}
\usepackage{amsmath}%
\providecommand{\U}[1]{\protect\rule{.1in}{.1in}}
\newtheorem{theorem}{Theorem}
\newtheorem{acknowledgement}[theorem]{Acknowledgement}

\newtheorem{corollary}[theorem]{Corollary}

\newtheorem{lemma}[theorem]{Lemma}

\newtheorem{proposition}[theorem]{Proposition}
\newtheorem{remark}[theorem]{Remark}

\begin{document}

\begin{center}
{\Large Convergence of series of conditional expectations}

\bigskip

M. Peligrad and C. Peligrad

\bigskip
\end{center}

Department of Mathematical Sciences

University of Cincinnati,

POBox 210025, Cincinnati, Oh 45221-0025, USA. \texttt{ }

E-mail address: peligrm@ucmail.uc.edu, peligrc@ucmail.uc.edu

\bigskip

\bigskip

\textbf{Abstract}

\bigskip

This paper deals with almost sure convergence for partial sums of Banach space
valued random variables. The results are then applied to solve similar
problems for weighted partial sums of conditional expectations. They are
further used to treat partial sums of powers of a reversible Markov chain
operator. The method of proof is based on martingale approximation. The
conditions are expressed in terms of moments of individual summands.

\bigskip

\noindent\textit{Keywords:} Almost sure convergence; Markov chains;
nonstationary sequences; maximal inequalities; smooth Banach spaces.

\smallskip

\noindent\textit{Mathematical Subject Classification} (2020): 60F15, 60E15, 60J05

\section{Introduction}

Let $\left(  X_{j}\right)  _{j\geq1}$ be a sequence random variables defined
on a probability space $(\Omega,\mathcal{F},P)$ with values in a separable
Banach space, adapted to $(\mathcal{F}^{j})_{j\geq1}$, a decreasing sequence
of sub-sigma algebras of $\mathcal{F}.$ We are going to study the convergence
of series $S_{n}=\sum_{j=1}^{n}X_{j}$. For integrable $X$, define{\LARGE \ }%
the conditional expectation $E^{j}X=E(X|\mathcal{F}^{j})$. Special attention
will be given to the situation when $X_{j}=a_{j}E(X|\mathcal{F}^{j})$ for
$\left(  a_{j}\right)  _{j\geq1}$ a sequence of real constants. Given a
stationary and reversible Markov chain $(\xi_{i})_{i\in Z}$ with values in a
measurable space $(S,\mathcal{S}),$ for an integrable function $f$ defined on
$S$ with values in a separable Banach space, let $X_{j}=f(\xi_{j})$. We denote
by $Q^{k}f(\xi_{0})=E(X_{k}|\xi_{0}),$ and also derive similar results for
$\sum_{j=1}^{k}a_{j}Q^{j}f$. It should be noted that Cohen, Cuny and Lin
(2017) studied the almost sure convergence of $L_{p}$ contractions for the
series of the form $\sum_{k=1}^{n}a_{k}Q^{k}f$ where $a_{n}$ is a Kaluza
sequence with divergent sum, $Q$ a power bounded operator and $\sum_{k\geq
1}\beta_{k}z^{k}$ converges in the open unit disk. We consider a general
sequence of constants and impose our conditions on the moments of $|E^{k}X|.$

The motivation for this study comes from a remarkable result, Theorem 3.11 in
Derriennic and Lin (2001). In the context of additive functionals of
stationary reversible Markov chains, for $f$ centered at expectation and
square integrable, if $E(S_{n}^{2})/n\rightarrow\sigma^{2}$ then we have that
$\sum_{k=1}^{n}k^{-1/2}Q^{k}f\ $converges a.s. Theorem 3.9 in the same paper
has a similar result for Harris recurrent Markov operators and some special
functions $f$.

As in Dedecker and Merlev\`{e}de (2008)\ or Cuny (2015), whenever
possible,\ we shall work with variables in a separable Banach space $B$. For
$x\in B$ for simplicity we shall denote the norm by $|x|=|x|_{B}.$

We denote by $L_{p}$ the set of measurable functions $X$ defined on a
probability space, with values in a separable Banach space such that
$||X||_{p}^{p}=E|X|^{p}<\infty$. For random variables $X$ in $L_{p}\,$the
notation $||X||_{p}\ll b$ means that there is a constant $C_{p}$ such that
$||X||_{p}\leq C_{p}b.$ For a sequence of positive constants $(a_{n})_{n},$
$(b_{n})_{n\text{ }},$ $a_{n}\ll b_{n}$ means that there is a constant $C$
such that $a_{n}\leq Cb_{n}.$

For integrable $X$, recall the notation{\LARGE \ }$E^{j}X=E(X|\mathcal{F}%
^{j}),$ $j\geq1$. For these definitions we direct the reader to the book of
Ledoux and Talagrand (1991). Now denote the reverse martingale difference
adapted to $(\mathcal{F}^{i})_{i\geq1}$ by
\begin{equation}
P^{i}(X)=E^{i}X-E^{i+1}X. \label{projdecr}%
\end{equation}

Sometimes we shall assume in addition that the Banach space is separable and
$r$-smooth for an $r$ such that $1<r\leq2.$ We shall use this property or
rather its consequence, for any sequence of $B$-valued martingale differences
$(X_{i})_{i\geq1},$ if $B$ is separable and $r$-smooth then for some $D>0,$%
\begin{equation}
E|X_{1}+X_{2}+...+X_{n}|^{r}\leq D(E|X_{1}|^{r}+E|X_{2}|^{r}+...+E|X_{n}%
|^{r}). \label{r-smooth}%
\end{equation}
(see Assouad, 1975).

We should mention that, according to our knowledge, our results are also new
for real-valued random variables.

The paper is divided in three parts. In Section 2 we obtain new maximal
inequalities. Then, in Section 3, we apply these maximal inequalities to
obtain convergence of series. Finally, in Section 4,  we apply the results to
reversible Markov chains.

\section{Maximal inequalities}

\begin{proposition}
\label{prop max ineq}Let $p>1$ and let $(X_{j})_{j\geq1}$ be a sequence of
Banach space valued random variables in $L_{p}$, adapted to a sequence of
decreasing sub-sigma fields of $\mathcal{F}$, $(\mathcal{F}^{k})_{k\geq1}%
.\ $Then
\[
E\max_{1\leq k\leq n}\left\vert S_{k}\right\vert ^{p}\ll E\max_{1\leq k\leq
n}|E^{k}S_{k}|^{p}+E\left\vert S_{n}\right\vert ^{p}.
\]
For $1\leq p\leq2$ and if $B$ is separable and $p$-smooth, we also have
\[
E\max_{1\leq k\leq n}\left\vert S_{k}\right\vert ^{p}\ll E\max_{1\leq k\leq
n}|E^{k}S_{k}|^{p}+\sum\nolimits_{i=1}^{n-1}E|P^{i}(S_{i})|^{p\text{ }}.
\]

\end{proposition}

\textbf{Proof}

Assume $S_{0}=0.$ It is easy to see that
\begin{equation}
S_{n}=E^{n}S_{n}+\sum\nolimits_{i=1}^{n-1}\left(  E^{i}S_{i}-E^{i+1}%
S_{i}\right)  . \label{cond exp decr}%
\end{equation}
By taking the maximum%

\begin{equation}
\max_{1\leq k\leq n}\left\vert S_{k}\right\vert \leq\max_{1\leq k\leq n}%
|E^{k}S_{k}|+\max_{2\leq k\leq n}\left\vert \sum\nolimits_{i=1}^{k-1}%
P^{i}(S_{i})\right\vert , \label{maxrep}%
\end{equation}
whence, by Minkowski's type inequality for norms in $L_{p}$ we obtain
\[
E\max_{1\leq k\leq n}\left\vert S_{k}\right\vert ^{p}\ll E\max_{1\leq k\leq
n}|E^{k}S_{k}|^{p}+E\max_{2\leq k\leq n}\left\vert \sum\nolimits_{i=1}%
^{k-1}P^{i}(S_{i})\right\vert ^{p}.
\]
By writing%
\[
\left\vert \sum\nolimits_{i=1}^{k-1}P^{i}(S_{i})\right\vert ^{p}\ll\left\vert
\sum\nolimits_{i=1}^{n-1}P^{i}(S_{i})\right\vert ^{p}+\left\vert
\sum\nolimits_{i=k}^{n-1}P^{i}(S_{i})\right\vert ^{p},
\]
we deduce that%
\begin{align}
E\max_{2\leq k\leq n}\left\vert \sum\nolimits_{i=1}^{k-1}P^{i}(S_{i}%
)\right\vert ^{p}  &  \ll E\left\vert \sum\nolimits_{i=1}^{n-1}P^{i}%
(S_{i})\right\vert ^{p}\label{first max}\\
&  +E\max_{1\leq k\leq n-1}\left\vert \sum\nolimits_{i=k}^{n-1}P^{i}%
(S_{i})\right\vert ^{p}.\nonumber
\end{align}
Since $P^{k}(S_{k})=E^{k}S_{k}-E^{k+1}S_{k}$ is a reverse martingale
difference adapted to the decreasing sequence of sigma algebras $(\mathcal{F}%
^{k})_{k\geq1}$, by Doob's maximal inequality for submartingales, for $p>1$%
\begin{equation}
E\max_{1\leq k\leq n-1}\left\vert \sum\nolimits_{i=k}^{n-1}P^{i}%
(S_{i})\right\vert ^{p}\ll E\left\vert \sum\nolimits_{i=1}^{n-1}P^{i}%
(S_{i})\right\vert ^{p}. \label{Doob}%
\end{equation}
By relation (\ref{cond exp decr}),%
\[
E\left\vert \sum\nolimits_{i=1}^{n-1}P^{i}(S_{i})\right\vert ^{p}\ll
E|E^{n}S_{n}|^{p}+E\left\vert S_{n}\right\vert ^{p},
\]
and the first inequality in this proposition follows.

The second inequality follows from (\ref{maxrep}) and (\ref{Doob}) combined
with (\ref{r-smooth}).

$\square$

\begin{remark}
\label{rem square}Note that if $p=2$ and if the variables have values in a
separable Hilbert space, then, for all $n\geq1$,
\[
E\left\vert \sum\nolimits_{i=1}^{n}P^{i}(S_{i})\right\vert ^{2}=\sum
\nolimits_{i=1}^{n}\left(  E|E^{i}S_{i}|^{2}-E|E^{i+1}S_{i}|^{2}\right)  .
\]

\end{remark}

We treat next linear combinations $%
{\textstyle\sum\nolimits_{j=1}^{k}}
a_{j}E^{j}X$ with $a_{j}$ a sequence of constants. We take $a_{j}$ real, but
complex valued constants can be treated in the same way. We shall assume that
$%
{\textstyle\sum\nolimits_{i=1}^{\infty}}
|a_{i}|=\infty$ since otherwise, by Doob's maximal inequality, we immediately
get for $p>1$
\[
E\max_{1\leq k\leq n}|%
{\textstyle\sum\nolimits_{j=1}^{n}}
a_{j}E^{j}X|^{p}\ll E\max_{1\leq k\leq n}|E^{k}X|^{p}\leq E|X|^{p},
\]
which is finite as soon as $X\in L_{p}$.

We shall also use the following notations
\[
s_{j}=\sum\nolimits_{i=1}^{j}a_{i},\text{ }s_{0}=0,\text{ }s_{k}^{\ast}%
=\max_{1\leq j\leq k}|\sum\nolimits_{i=1}^{j}a_{i}|,
\]
and
\begin{equation}
b_{k}=\max\left(  k^{-1}\left(  s_{4k}^{\ast}\right)  ^{2},s_{k}^{2}%
-s_{k-1}^{2}\right)  . \label{defbk}%
\end{equation}
The next corollary follows from Proposition \ref{prop max ineq} applied to
$X_{i}=a_{i}E^{i}X$.

\begin{corollary}
\label{cor MAX}Let $X$ be a random variable with values in a separable Banach
space and, for any $p>1,$ $E|X|^{p}<\infty$. Then%
\begin{equation}
E\max_{1\leq k\leq n}\left\vert \sum\nolimits_{j=1}^{k}a_{j}E^{j}X\right\vert
^{p}\ll E\max_{1\leq k\leq n}|s_{k}E^{k}X|^{p}+E\left\vert \sum\nolimits_{j=1}%
^{n}a_{j}E^{j}X\right\vert ^{p}.\label{maxlin}%
\end{equation}
For $1\leq p\leq2$ and if $B$ is $p$-smooth, we also have
\begin{equation}
E\max_{1\leq k\leq n}\left\vert \sum\nolimits_{j=1}^{k}a_{j}E^{j}X\right\vert
^{p}\ll E\max_{1\leq k\leq n}|s_{k}E^{k}X|^{p}+\sum\nolimits_{i=1}%
^{n-1}E|s_{i}P^{i}(X)|^{p\text{ }}.\label{maxlin2}%
\end{equation}

\end{corollary}

\begin{remark}
\label{rem-dyadic} For $p\geq1$, an estimate of $E\max_{1\leq k\leq n}%
|s_{k}E^{k}X|^{p}$ is
\[
E\max_{1\leq k\leq n}|s_{k}E^{k}X|^{p}\ll\sum\nolimits_{k=1}^{n}k^{-1}\left(
s_{4k}^{\ast}\right)  ^{p}E|E^{k}X|^{p}.
\]
\textbf{ }
\end{remark}

\textbf{Proof of Remark \ref{rem-dyadic}}

\bigskip

By Doob's maximal inequality
\[
E\max_{2^{i}\leq k\leq2^{i+1}}|E^{k}X|^{p}\ll E|E^{2^{i}}X|^{p}.
\]
Note that
\begin{gather*}
E\max_{1\leq k\leq2^{r}}\left\vert s_{k}E^{k}X\right\vert ^{p}\leq
\sum\nolimits_{i=0}^{r-1}E\max_{2^{i}\leq k\leq2^{i+1}}|s_{k}E^{k}X|^{p}\\
\leq\sum\nolimits_{i=0}^{r-1}\left(  s_{2^{i+1}}^{\ast}\right)  ^{p}%
E\max_{2^{i}\leq k\leq2^{i+1}}|E^{k}X|^{p}\ll\sum\nolimits_{i=0}^{r-1}\left(
s_{2^{i+1}}^{\ast}\right)  ^{p}E|E^{2^{i}}X|^{p}.
\end{gather*}
Using the fact that $(s_{n}^{\ast})_{n\geq1}$ is increasing and $\left(
E|E^{n}X|^{p}\right)  _{n\geq1}$ is decreasing, we easily obtain%
\begin{gather*}
E\max_{1\leq k\leq2^{r}}\left\vert s_{k}E^{k}X)\right\vert ^{p}\ll
\sum\nolimits_{i=0}^{r-1}\left(  s_{2^{i+1}}^{\ast}\right)  ^{p}E|E^{2^{i}%
}X|^{p}\\
\ll\sum\nolimits_{i=1}^{2^{r-1}}i^{-1}\left(  s_{4i}^{\ast}\right)
^{p}E|E^{i}X|^{p}\ll\sum\nolimits_{i=1}^{2^{r-1}}i^{-1}\left(  s_{4i}^{\ast
}\right)  ^{p}E|E^{i}X|^{p}.
\end{gather*}

Now if $2^{r-1}<n\leq2^{r},$ clearly
\[
E\max_{1\leq k\leq n}\left\vert s_{k}E^{k}X\right\vert ^{p}\leq E\max_{1\leq
k\leq2^{r}}\left\vert s_{k}E^{k}X\right\vert ^{p}\ll\sum\nolimits_{k=1}%
^{n}k^{-1}\left(  s_{4k}^{\ast}\right)  ^{p}E|E^{k}X|^{p}.
\]

$\square$

\bigskip

For $p=2$ we also have the following result:

\begin{corollary}
\label{Cor second moments}Assume that $X$ has values in a separable Hilbert
space and $E|X|^{2}<\infty$. Then
\[
E\max_{1\leq k\leq n}\left\vert \sum\nolimits_{j=1}^{k}a_{j}E^{j}X\right\vert
^{2}\ll\sum\nolimits_{k=1}^{n}b_{k}E|E^{k}X|^{2}.
\]

\end{corollary}

\begin{remark}
\label{rem 2 particular}In particular, under the conditions of Corollary
\ref{Cor second moments}, if $a_{j}=j^{-1/2}$ for all $j\geq1,$
\[
E\max_{1\leq k\leq n}\left\vert \sum\nolimits_{j=1}^{k}j^{-1/2}E^{j}%
X\right\vert ^{2}\ll\sum\nolimits_{j=1}^{n}E|E^{j}X|^{2},
\]
and if $a_{j}=1$ for all $j\geq1$, then
\[
E\max_{1\leq k\leq n}\left\vert \sum\nolimits_{j=1}^{k}E^{j}X\right\vert
^{2}\ll\sum\nolimits_{j=1}^{n}jE|E^{j}X|^{2}.
\]

\end{remark}

\textbf{Proof of Corollary \ref{Cor second moments}}.

\bigskip

In a Hilbert space $P^{i}(X)$ and $P^{j}(X)$ are orthogonal for $i\neq j$. By
the properties of the conditional expectations and Remark \ref{rem square},
$(s_{0}=0,$ $n\geq2)$%
\begin{gather*}
E\left\vert \sum\nolimits_{i=1}^{n-1}s_{i}P^{i}(X)\right\vert ^{2}%
=\sum\nolimits_{i=1}^{n-1}s_{i}^{2}E|P^{i}(X)|^{2}\\
=\sum\nolimits_{i=1}^{n-1}s_{i}^{2}\left(  E|E^{i}X|^{2}-E|E^{i+1}%
X|^{2}\right)  \leq\sum\nolimits_{j=1}^{n-1}\left(  s_{j}^{2}-s_{j-1}%
^{2}\right)  E|E^{j}X|^{2}.
\end{gather*}
Combining the latter inequality with Proposition \ref{prop max ineq} and
Remark \ref{rem-dyadic}, we get
\begin{align*}
E\max_{1\leq k\leq n}\left\vert \sum\nolimits_{j=1}^{k}a_{j}E^{j}X\right\vert
^{2}  &  \ll\sum\nolimits_{k=1}^{n}k^{-1}\left(  s_{4k}^{\ast}\right)
^{2}E|E^{k}X|^{2}\\
&  +\sum\nolimits_{k=1}^{n-1}\left(  s_{k}^{2}-s_{k-1}^{2}\right)
E|E^{k}X|^{2}\\
&  \leq\sum\nolimits_{k=1}^{n}b_{k}E|E^{k}X|^{2}.
\end{align*}
where $b_{k}$ is given in (\ref{defbk}).

$\square$

\section{Convergence of series}

We give some straightforward applications of the maximal inequalities
established before.

\begin{proposition}
\label{propconv}Let $p>1$ and let $(X_{j})_{j\geq1}$ be a sequence of Banach
space valued random variables. Assume that $(S_{n})_{n\geq1}$ is bounded in
$L_{p}$ and $\left(  E^{n}S_{n}\right)  _{n\geq1}$ converges in $L_{p}$. Then
$\left(  S_{n}\right)  _{n\geq1}$ converges in $L_{p}$. If in addition
$\left(  E^{n}X_{n}\right)  _{n\geq1}$ converges a.s. then $\left(
S_{n}\right)  _{n\geq1}$ converges $a.s.$
\end{proposition}

\textbf{Proof of Proposition \ref{propconv}}.

\bigskip

We start from the representation given in (\ref{cond exp decr}).

Now note that $\sum\nolimits_{i=1}^{n-1}\left(  E^{i}S_{i}-E^{i+1}%
S_{i}\right)  $ is a reversed martingale, which converges a.s. and in $L_{p}$ provided%

\[
\sup_{n}E\left\vert \sum\nolimits_{i=1}^{n-1}\left(  E^{i}S_{i}-E^{i+1}%
S_{i}\right)  \right\vert ^{p}<\infty.
\]
By (\ref{cond exp decr}) and the triangle inequality this condition is
satisfied under the conditions of this proposition, and the result follows.

$\square$

\bigskip

If we apply the previous proposition to $X_{i}=a_{i}X$ we obtain the following
corollary. We denote as before $s_{n}=\sum\nolimits_{i=1}^{n}a_{i}.$

\begin{corollary}
\label{a.s. conv series p}Let $p>1$. Assume that $(%
{\textstyle\sum\nolimits_{i=1}^{n}}
a_{i}E^{i}X)_{n\geq1}$ is bounded in $L_{p}$ and $\left(  s_{n}E^{n}X\right)
_{n\geq1}$ converges in $L_{p}.$ Then $\left(
{\textstyle\sum\nolimits_{i=1}^{n}}
a_{i}E^{i}X\right)  _{n\geq1}$ converges in $L_{p}$. If in addition $\left(
s_{n}E^{n}X\right)  _{n\geq1}$ converges a.s. then $\left(
{\textstyle\sum\nolimits_{i=1}^{n}}
a_{i}E^{i}X\right)  _{n\geq1}$ converges $a.s.$
\end{corollary}

\begin{corollary}
\label{a.s. conv series 2}Assume that $X$ has values in a separable Hilbert
space. Define $(b_{k})$ by (\ref{defbk}), and assume that
\[
\sum\nolimits_{k\geq1}b_{k}E|E^{k}X|^{2}<\infty.
\]
Then $\left(
{\textstyle\sum\nolimits_{i=1}^{n}}
a_{i}E^{i}X\right)  _{n\geq1}$ converges in $L_{2}$ and $a.s.$
\end{corollary}

\textbf{Proof of Corollary \ref{a.s. conv series 2}}. We start from the
representation (\ref{cond exp decr}), namely
\[
\sum\nolimits_{j=1}^{k}a_{j}E^{j}X=s_{k}E^{k}X+\sum\nolimits_{i=1}^{k-1}%
s_{i}\left(  E^{i}X-E^{i+1}X\right)  .
\]
The reverse martingale $\sum\nolimits_{i=1}^{k-1}s_{i}\left(  E^{i}%
X-E^{i+1}X\right)  $ converges a.s. and in $L_{2}$ provided is bounded in
$L_{2}.$ Because the variables have values in a Hilbert space, by Remark
\ref{rem square}%
\[
E\left\vert \sum\nolimits_{i=1}^{k-1}s_{i}\left(  E^{i}X-E^{i+1}X\right)
\right\vert ^{2}=\sum\nolimits_{i=1}^{k-1}s_{i}^{2}\left(  E|E^{i}%
X|^{2}-E|E^{i+1}X|^{2}\right)  ,
\]
which is bounded because  $\sum\nolimits_{k=1}^{n-1}\left(  s_{k}^{2}%
-s_{k-1}^{2}\right)  E|E^{k}X|^{2}$ is positive and we assumed it is bounded.

By a similar proof as of Corollary \ref{Cor second moments} because
$(s_{k}^{\ast})_{k\text{ }}$is increasing and $(E^{k}X)_{k}$ is decreasing, we
obtain
\begin{gather*}
E(\max_{2^{a}\leq k\leq2^{b}}|s_{k}E^{k}X|^{2})\leq\sum\nolimits_{j=a}%
^{b-1}E(\max_{2^{i}\leq k\leq2^{i+1}}|s_{k}E^{k}X|^{2})\\
\ll\sum\nolimits_{j=a}^{b-1}\left(  s_{2^{j+1}}^{\ast}\right)  ^{2}E|E^{2^{j}%
}X|^{2}\ll\sum\nolimits_{j=2^{a-1}}^{2^{b}}j^{-1}\left(  s_{4j}^{\ast}\right)
^{2}E|E^{j}X|^{2}.
\end{gather*}
Now if $2^{a}\leq m<2^{a+1}$ and $2^{b-1}\leq n<2^{b},$ with $a$, $b$ integers
$a\leq b,$ then%
\[
E(\max_{m\leq k\leq n}|s_{k}E^{k}X|^{2})\ll\sum\nolimits_{j=2^{a}}^{2^{b}%
}j^{-1}\left(  s_{4j}^{\ast}\right)  ^{2}E|E^{j}X|^{2}.
\]
This implies relation (22.11) in Billingsley (1999) and the proof continues as there.

$\square$

\section{\textbf{Reversible Markov chains.}}

All the results presented so far, are useful to treat power series of
operators associated to stationary reversible Markov chains. Let $(\xi
_{i})_{i\in Z}$ be defined on $(\Omega,\mathcal{F},P)$ with values in a
general measurable space $(S,\mathcal{S},\pi)$ with stationary transitions
$Q(x,A)$. Reversible, means that the distribution of $(\xi_{i},\xi_{i+1})$ is
the same as the distribution of $(\xi_{i+1},\xi_{i}).$ For $f$ defined on $S$
with values in a separable Banach space $B$, and any $n\in N$ denote by
$Q^{n}f(\xi_{0})=E(f(\xi_{n})|\xi_{0}).$

We also denote by $Q$ the operator on integrable $f$ defined by
\[
Qf(x)=\int\nolimits_{S}f(y)Q(x,dy).
\]
We denote the invariant distribution by $\pi,$ which is a measure on
$\mathcal{S}.$ The integral with respect to $\pi$ is denoted by $E_{\pi}.$

We shall use below notations similar to those in previous sections,
\begin{align*}
s_{k}^{e}  &  =a_{2}+...+a_{2k},\text{ }s_{k}^{e\ast}=\max_{1\leq j\leq
k}|s_{j}^{e}|\\
b_{k}^{e}  &  =\max\left(  k^{-1}\left(  s_{4k}^{e\ast}\right)  ^{2},\left(
s_{k}^{e}\right)  ^{2}-\left(  s_{k-1}^{e}\right)  ^{2}\right)
\end{align*}%
\begin{align*}
s_{k}^{0}  &  =a_{1}+...+a_{2k-1},s_{k}^{o\ast}=\max_{1\leq j\leq k}|s_{j}%
^{o}|\\
\text{ }b_{k}^{0}  &  =\max\left(  k^{-1}\left(  s_{4k}^{o\ast}\right)
^{2},\left(  s_{k}^{o}\right)  ^{2}-\left(  s_{k-1}^{o}\right)  ^{2}\right)
\end{align*}%
\[
b_{k}^{\ast}=\max\left(  b_{k}^{e},b_{k}^{o}\right)  .
\]

\begin{theorem}
\label{Th max ineq rev two}For $f$ with values in a separable Hilbert space
with $E_{\pi}|f|^{2}<\infty$, we have
\[
E_{\pi}\max_{1\leq k\leq2n}\left\vert \sum\nolimits_{j=1}^{k}a_{j}%
Q^{j}f\right\vert ^{2}\ll\sum\nolimits_{j=1}^{n}b_{j}^{\ast}E_{\pi}%
|Q^{j}f|^{2}.
\]

\end{theorem}

\begin{corollary}
In particular%
\[
E_{\pi}\max_{1\leq k\leq n}\left\vert \sum\nolimits_{j=1}^{k}Q^{j}f\right\vert
^{2}\ll\sum\nolimits_{j=1}^{n}jE_{\pi}|Q^{j}f|^{2}%
\]
and
\[
E_{\pi}\max_{1\leq k\leq n}\left\vert \sum\nolimits_{j=1}^{k}j^{1/2}%
Q^{j}f\right\vert ^{2}\ll\sum\nolimits_{j=1}^{n}E_{\pi}|Q^{j}f|^{2}.
\]

\end{corollary}

To relate our result with Theorem 3.11 in Derriennic and Lin (2001), we give
the following result similar to their theorem.

\begin{corollary}
\label{corineqreal}If $f$ is real valued, mean $0$ and has finite second
moment, then%
\[
E_{\pi}\sup_{k\geq1}\left\vert \sum\nolimits_{j=1}^{k}j^{1/2}Q^{j}f\right\vert
^{2}<\infty
\]
provided one of the following equivalent conditions hold
\begin{equation}
(a)\text{ \ \ }\sum\nolimits_{k=1}^{n}E_{\pi}(fQ^{k}f)\text{ is bounded in
}L_{2}\text{.}\label{condrevreal}%
\end{equation}%
\[
(b)\text{\ \ \ \ \ \ \ \ \ \ \ \ \ \ \ \ \ \ \ \ \ \ \ \ \ }\sup_{n}ES_{n}%
^{2}/n<C.
\]%
\[
(c)\ \ \ \ \ \ \ \ \ \ \ \ \ \ \ \ \ \ \ \ \ \ \ \lim_{n\rightarrow\infty
}ES_{n}^{2}/n=\sigma^{2}.
\]%
\[
(d)\text{ \ \ \ \ \ \ \ \ \ \ \ \ \ \ \ \ \ \ \ \ \ \ \ \ }\int\limits_{-1}%
^{1}\frac{1}{1-t}d\rho_{f}<\infty.
\]%
\[
(e)\text{ \ \ \ \ \ \ \ \ \ \ \ \ \ \ \ \ \ \ \ \ \ \ \ \ \ \ \ \ }f\in
\sqrt{1-Q}L_{2}^{0}.
\]
Above, $\rho_{f}$ denotes the spectral measure of $f$ associated with the
self-adjoint operator $Q$ and function $f$ on $L_{2}(S,\pi)$. Also $L_{2}^{0}$
is the set of functions which are square integrable and have mean $0.$
\end{corollary}

The proofs of the results in this section are based on the following identity,
relating $Q^{j}f$ to $E^{j}f.$

\begin{lemma}
\label{rel with reversible}For $n\geq1$%
\[
\sum\nolimits_{j=1}^{2n}a_{j}Q^{j}f(\xi_{0})=E_{0}\left(  \sum\nolimits_{j=1}%
^{n}a_{2j}E^{j}f\left(  \xi_{0}\right)  +\sum\nolimits_{j=0}^{n-1}%
a_{2j+1}E^{j+1}f(\xi_{1})\right)  .
\]

\end{lemma}

\textbf{Proof of Lemma \ref{rel with reversible}}. We estimate $E(X_{2j}%
|\mathcal{F}_{0})=Q^{2j}f(\xi_{0}).$ By the Markov property and reversibility%
\begin{align*}
a_{2j}Q^{2j}f(\xi_{0})  &  =a_{2j}E_{0}E(f(\xi_{2j})|\xi_{j})\\
&  =a_{2j}E_{0}E(f(\xi_{0})|\mathcal{F}^{j})=a_{2j}E_{0}E^{j}f(\xi_{0}).
\end{align*}
By the above identity
\begin{equation}
\sum\nolimits_{j=1}^{n}a_{2j}Q^{2j}f(\xi_{0})=E_{0}\sum\nolimits_{j=1}%
^{n}a_{2j}E^{j}f(\xi_{0}). \label{even}%
\end{equation}
Similarly%
\begin{gather*}
a_{2j+1}(Q^{2j+1}f)(\xi_{0})=a_{2j+1}E(f(\xi_{2j+1})|\xi_{0})=\\
a_{2j+1}E_{0}E(f(\xi_{2j+1})|\xi_{j+1})=a_{2j+1}E_{0}E(f(\xi_{1})|\xi
_{j+1})=a_{2j+1}E_{0}E^{j+1}f(\xi_{1}),
\end{gather*}
and so
\begin{equation}
\sum\nolimits_{j=0}^{n-1}a_{2j+1}Q^{2j+1}f(\xi_{0})=E_{0}\sum\nolimits_{j=0}%
^{n-1}a_{2j+1}E^{j+1}f(\xi_{1}). \label{odd}%
\end{equation}
Overall, by (\ref{even}) and (\ref{odd}), we have the result of this lemma.

$\square$

\begin{remark}
\label{rem rel with reversible}A similar result as in Lemma
\ref{rel with reversible} holds for odd sums. We easily deduce%
\begin{align*}
\max_{1\leq k\leq2n}\left\vert \sum\nolimits_{j=1}^{k}a_{j}Q^{j}f(\xi
_{0})\right\vert  &  \leq E_{0}\max_{1\leq k\leq n}\left\vert \sum
\nolimits_{j=1}^{k}a_{\ j}E^{j}f(\xi_{0})\right\vert \\
&  +E_{0}\max_{1\leq k\leq n}\left\vert \sum\nolimits_{j=1}^{k}a_{2j+1}%
E^{j}f(\xi_{1})\right\vert .
\end{align*}
whence, for every $p\geq1,$
\begin{align*}
E\max_{1\leq k\leq2n}\left\vert \sum\nolimits_{j=1}^{k}a_{j}Q^{j}f(\xi
_{0})\right\vert ^{p}  &  \ll E\max_{1\leq k\leq n}\left\vert \sum
\nolimits_{j=1}^{k}a_{2j}E^{j}f(\xi_{0})\right\vert ^{p}\\
&  +E\max_{1\leq k\leq n}\left\vert \sum\nolimits_{j=1}^{k}a_{2j+1}E^{j}%
f(\xi_{0})\right\vert ^{p}.
\end{align*}

\end{remark}

It is straightforward now to combine the Remark \ref{rem rel with reversible}
with the results in the previous sections. It can be easily combined to
Corollary \ref{cor MAX}.

In particular, for $p=2$ we combine Remark \ref{rem rel with reversible} with
Corollary \ref{Cor second moments} and we get the result in Theorem
\ref{Th max ineq rev two}.

\bigskip

\textbf{Proof of Corollary \ref{corineqreal}}

\bigskip

First we find a more flexible maximal inequality which has interest in itself.
From Theorem \ref{Th max ineq rev two} with $f$ replaced by $f+Qf,$ we get%

\[
E_{\pi}\max_{1\leq k\leq2n}\left(  \sum\nolimits_{j=1}^{k}(Q^{j}%
f+Q^{j+1}f)\right)  ^{2}\ll\sum\nolimits_{j=1}^{n}jE_{\pi}\left(
Q^{j}f+Q^{j+1}f\right)  ^{2}.
\]
But, by the fact that on $L_{2}$ the operator $Q$ is self-adjoint, with the
notation $<f,g>=E_{\pi}fg,$
\[
E_{\pi}\left(  Q^{j}f+Q^{j+1}f\right)  ^{2}=<f,Q^{2j}f>+2<f,Q^{2j+1}%
>+<f,Q^{2j+2}f>.
\]
So, after some algebraic computation, using that $<f,Q^{2j}f>=||Q^{j}%
f||^{2}\geq0$ we obtain%
\[
\sum\nolimits_{j=1}^{k}jE_{\pi}\left(  Q^{j}f+Q^{j+1}f\right)  ^{2}%
\ll\left\vert \sum\nolimits_{j=1}^{2k}jE_{\pi}\left(  fQ^{j}f\right)
\right\vert +<f,Q^{2}f>,
\]
and therefore%
\begin{equation}
E_{\pi}\max_{1\leq k\leq2n}\left(  \sum\nolimits_{j=1}^{k}(Q^{j}%
f+Q^{j+1}f)\right)  ^{2}\ll\left\vert \sum\nolimits_{j=1}^{2n}jE_{\pi}\left(
fQ^{j}f\right)  \right\vert +<f,Q^{2}f>. \label{est partial}%
\end{equation}
On the other hand, because
\[
\sum\nolimits_{j=1}^{k}(Q^{j}f+Q^{j+1}f)=2\sum\nolimits_{j=1}^{k}%
Q^{j}f-Qf+Q^{k+1}f,
\]
we obtain
\[
\max_{1\leq k\leq2n}\left\vert \sum\nolimits_{j=1}^{k}Q^{j}f\right\vert
\ll\max_{1\leq k\leq2n}\left\vert \sum\nolimits_{j=1}^{k}(Q^{j}f+Q^{j+1}%
f)\right\vert +\max_{1\leq k\leq2n}|Q^{k+1}f|+|Qf|.
\]
To estimate the second moment of $\max_{1\leq k\leq n+1}|Q^{k}f|$ we use the
Stein Theorem (see page 106 in Stein (1970) or Krengel (1895), page 190).
\[
E_{\pi}\max_{1\leq k\leq2n}|Q^{k+1}f|^{2}\leq E_{\pi}(Qf)^{2}=E_{\pi}%
(fQ^{2}f).
\]
By combining these results with the estimate in (\ref{est partial}) we get
\[
E_{\pi}\max_{1\leq k\leq2n}\left\vert \sum\nolimits_{j=1}^{k}Q^{j}f\right\vert
^{2}\ll\left\vert \sum\nolimits_{j=1}^{2n}jE_{\pi}\left(  fQ^{j}f\right)
\right\vert +E_{\pi}(fQ^{2}f),
\]
and the result follows.

The equivalences (a)-(e) are well-known results in the literature (see for
instance pages 3 and 4 in Kipnis and Varadhan (1986) and Cuny (2009) ).

$\square$

\bigskip\ 

\begin{acknowledgement}
This paper was partially supported by the NSF grant DMS-2054598. The authors
would like to thank Christophe Cuny for point out a gap in a former version of
the paper.
\end{acknowledgement}

\textbf{\ }
\end{document}